\documentclass[pdflatex,sn-mathphys-num]{sn-jnl}

\usepackage{graphicx}%
\usepackage{multirow}%
\usepackage{amsmath,amssymb,amsfonts}%
\usepackage{amsthm}%
\usepackage{mathrsfs}%
\usepackage[title]{appendix}%
\usepackage{xcolor}%
\usepackage{textcomp}%
\usepackage{manyfoot}%
\usepackage{booktabs}%
\usepackage{algorithm}%
\usepackage{algorithmicx}%
\usepackage{algpseudocode}%
\usepackage{listings}%

\usepackage{mathtools}%

\theoremstyle{plain}
 \newtheorem{thm}{Theorem}[section]
 \newtheorem{prop}[thm]{Proposition}
 
 \newtheorem{cor}[thm]{Corollary}

\theoremstyle{definition}

 \newtheorem*{dfn*}{Definition}

 \newtheorem*{hyp*}{Hypothesis}
 \newtheorem*{hyps*}{Hypotheses}
 
\theoremstyle{remark}
 \newtheorem*{rem*}{Remark}
 \newtheorem*{warn*}{Warning}
 \newtheorem*{rems*}{Remarks}
 \newtheorem*{not*}{Notation}
 \newtheorem*{nots*}{Notations}
 \newtheorem*{convs*}{Conventions}
 \newtheorem*{exm*}{Example}
 \newtheorem*{pf*}{Proof}
 \numberwithin{equation}{section}

\renewcommand{\leq}{\leqslant}
\renewcommand{\geq}{\geqslant}

\newcommand{\QQ}{\mathbb{Q}}
\newcommand{\CC}{\mathbb{C}}
\newcommand{\ZZ}{\mathbb{Z}}
\newcommand{\NN}{\mathbb{N}}

\newcommand{\Cp}{\mathbb{C}_p}
\newcommand{\Zp}{\mathbb{Z}_p}

\newcommand{\Bf}{\mathbb{F}}

\newcommand{\ord}{\text{\rm ord}}

\begin{document}

\title[Special values of $p$-adic $L$-functions and $\lambda$-invariants]{Special values of $p$-adic $L$-functions and Iwasawa $\lambda$-invariants of Dirichlet characters}


\author*{\fnm{Heiko} \sur {Knospe}} 
\affil{\orgname{Technische Hochschule K\"{o}ln}, \orgdiv{Institute of Computer and Communication Technology}, \orgaddress{\street{Betzdorfer Str. 2}, \city{K\"{o}ln}, \postcode{50679}, \country{Germany}}}
\email{heiko.knospe@th-koeln.de}

\abstract{
We investigate the Iwasawa $\lambda$-invariant $\lambda_p(\chi)$ of the $p$-adic $L$-function associated with a Dirichlet character $\chi$ for odd primes $p$. Using special values of the $p$-adic $L$-function and its derivative, we derive novel and computationally efficient criteria to distinguish between the cases
$\lambda = 0$, $\lambda = 1$, $\lambda = 2$, and $\lambda \geq 3$.  In particular, we  study 
the case when the $p$-adic $L$-function vanishes at $s=0$.
Building on the results of Ferrero-Greenberg and of Gross-Koblitz,
we establish explicit conditions for $\lambda_p(\chi) >1 $ and
$\lambda_p(\chi) > 2$. Furthermore, we extend the methods of Ernvall-Metsänkylä and of Dummit et al.\ to calculate the $\lambda$-invariant. This is achieved by twisting $\chi$ by characters $\psi$ of $p$-power order and using the values of the $p$-adic $L$-function at $s=2-p, \dots, 0$. 
Additionally, we utilize the value at $s=1$ to compute $\lambda_p(\chi)$.
Finally, these results are applied to generate numerical data on the distribution of $\lambda$-invariants, when either the prime $p$ or the Dirichlet character $\chi$ is fixed. 
}


\maketitle

\section{Introduction}\label{Sect1}
Let $p$ be an odd prime and  let $\chi$ be a Dirichlet character.
The Kubota-Leopoldt $L$-function $L_p(s,\chi)$ is a meromorphic function of the $p$-adic variable $s$, which interpolates
special values of the complex $L$-function $L(s,\chi)$ at non-positive integers  \cite{KuLe}. Specifically, the interpolation formula is given by
\begin{equation}\label{InterpFormula}
L_p(1-n, \chi) \;=\; - \big( 1 -\chi \omega^{-n} (p) p^{n-1} \big) \times \frac{B_{n, \chi \omega^{-n} }}{n}
\quad\text{at every integer $n \geq 1$, }
\end{equation}
where $\omega$ denotes the Teichm\"{u}ller character modulo $p$ and $B_{n, \chi \omega^{-n}}$
 the $n$-th $\chi \omega^{-n}$-twisted Bernoulli number.

 In the following, we assume that $\chi$ is a {\em non-trivial even character of the first kind}, expressed as $\chi = \theta \omega^i$,
where the conductor of the character $\theta$ is not divisible by $p$ and $i \in \{0,\dots, p-2\}$. 
We fix an embedding of $\overline{\QQ}$ into $\Cp$.
We let $\zeta_n$ denote a primitive $n$-th root of unity.
We set $\mathcal{O}_{\chi}:=\Zp[\text{\rm Im}(\chi)] = \Zp[\zeta_{\ord(\chi)}]$, and let $f$ denote the residue class degree of $\mathcal{O}_{\chi}$.
Let $\psi$ be a {\em character of the second kind} which is either trivial or has $p$-power order. Then $\zeta_{\psi}= \psi(1+p)^{-1}$ is either $1$ or a root of unity of $p$-power order. 
Iwasawa associated a unique power series $\mathcal{F}_{\chi}(T) = \sum_{j=0}^{\infty} c_j T^j \in \mathcal{O}_{\chi}[\![T]\!]$ satisfying
\begin{equation}
\mathcal{F}_{\chi}\big(\zeta_{\psi} (1+ p)^{s} - 1\big) =\;
L_p(s, \chi \psi)
\quad\text{for all $s\in\Cp$ with $\big|s\big|_p\! <p^{\frac{p-2}{p-1}}$.}
\label{powerseries}
\end{equation}
The $\mu$-invariant of $\mathcal{F}_{\chi}$ vanishes by the Ferrero-Washington theorem \cite{FeWa}, and $\mathcal{F}_{\chi}$ factorises into a product of an invertible power series
and a uniquely determined {\em distinguished polynomial}. The degree of this polynomial,  $\lambda=\lambda_p(\chi)$, equals the number of zeroes of $\mathcal{F}_{\chi}(T)$ on the open unit disk.

In our article \cite{DeKn}, we generalized a conjecture of Ellenberg-Jain-Venkatesh \cite{EJV} and predicted that the probability  
 that $\lambda_p(\chi)=r$
is approximately equal to
\begin{equation}
p^{-fr}\times\;\prod_{t>r}\;\big(1-p^{-ft}\big).  \label{pred-prob}
\end{equation}

In the present article, we give several new formulas for the $\lambda$-invariant. In particular, we investigate the case where the $p$-adic $L$-function and the associated power series have a zero  at $s=0$ and $T=0$, respectively (the rank one case).
By a result of Ferrero-Greenberg \cite{FeGr}, the order of $L_p(s,\chi)$ at $s=0$ is zero when $\chi \omega^{-1}(p) \neq 1$ and one when $\chi \omega^{-1}(p) = 1$.

The main results of this article are Theorem \ref{theorem-trivial} (for the rank one case), Theorem \ref{theorem-general} (for the rank zero case), and Theorem \ref{thm1}. These formulas allow for a convenient computation of the $\lambda$-invariant without requiring deep knowledge of $p$-adic $L$-functions and Iwasawa power series.
We also provide numerical evidence that the prediction  (\ref{pred-prob}) on the distribution of $\lambda$-invariants remains valid in the rank one case if we shift $\lambda$ by $1$.

\section{Analyzing $\lambda$-invariants via special values of $L_p(s,\chi )$}
\label{Sect2}

The following Proposition provides an iterative method for determining the $\lambda$-invariant.

\begin{prop} Let $K/\QQ_p$ be a finite extension with ring of integers $\mathcal{O}_K$ and let $v_p$ be the valuation of $K$ normalized such that $v_p(p)=1$. Let
$\mathcal{F}(T)  \in \mathcal{O}_K[[T]]$ be a power series with  $\mu_p(\mathcal{F})=0$, and choose $t_0 \in \mathcal{O}_K$ satisfying $v_p(t_0) > 0$. 
\begin{enumerate}
\item[(i)] $\lambda_p (\mathcal{F}) > 0$ if and only if $v_p(\mathcal{F}(t_0))>0$. 
\item[(ii)] Let  $\mathcal{F}(T)  = \sum_{j=0}^{\infty} a_j (T-t_0)^j$ be the Taylor expansion of $\mathcal{F}(T)$ at $t_0$
and let $p_n(T) = a_0 + a_1 (T-t_0) + \dots + a_n (T-t_0)^n$. Choose any $b \in \mathcal{O}_K$ such that $v_p(b-t_0)=1$.
Then for any $n \geq 1$:
$$ \lambda_p (\mathcal{F}) > n \Longleftrightarrow \lambda_p (\mathcal{F}) > n-1 \text{ and } v_p( \mathcal{F}(b) - p_{n-1}(b) ) > n  .$$
\end{enumerate}
The latter condition can be verified by computing  $\mathcal{F}(b)$ and $p_{n-1}(b)$ up to precision of $O(p^{n+1})$.
\label{lambdacomp}
\end{prop}

\begin{proof} The $\lambda$-invariant of $\mathcal{F}(T)$ can be computed using an expansion at any $t_0 \in \mathcal{O}_K$ satisfying $v_p(t_0) > 0$. 
Part (i) follows from the definition of the $\lambda$-invariant.  For part (ii), let $n \geq 1$. Assuming $\lambda_p(\mathcal{F}) > n-1$, which implies
$v_p(a_j)>0$ for $j \leq n-1$, the condition $\lambda_p(\mathcal{F})>n$ is equivalent to $v_p(a_n)>0$. We have
$$ \mathcal{F}(b) - p_{n-1}(b) = a_n (b-t_0)^n + a_{n+1} (b-t_0)^{n+1} + \dots $$
Since $v_p(b-t_0)=1$, we have $v_p((b-t_0)^n)=n$ and $v_p(a_n)>0$ if and only if $v_p ( a_n (b-t_0)^n) > n$. This is equivalent to $v_p( \mathcal{F}(b) - p_{n-1}(b) ) > n$, as the valuation of $\sum_{j=n+1}^{\infty} a_j (b-t_0)^j$ is at least $n+1$.
\end{proof}

\begin{cor} Let $\mathcal{F}_{\chi}(T)$ be the  Iwasawa power series associated with $\chi$, and let $t_0, \ b \in \mathcal{O}_{\chi}$ be as in Proposition \ref{lambdacomp}.
Consider the conditions:\\
\begin{enumerate}
\item[(1)]   $v_p(\mathcal{F}_{\chi}(t_0)) > 0$ 
\item[(2)] $v_p(\mathcal{F}_{\chi}(b)-\mathcal{F}_{\chi}(t_0)) > 1$
\item[(3)]  $v_p(\mathcal{F}_{\chi}(b)- \mathcal{F}_{\chi}'(t_0) (b-t_0) - \mathcal{F}_{\chi}(t_0)) > 2$
\end{enumerate}
Then (1) can be checked modulo $p$, (2) is equivalent to $v_p(\mathcal{F}'_{\chi}(t_0)) > 0$ and can be checked
modulo $p^2$, and (3) can be checked modulo $p^3$. Moreover,
\begin{align*}
\lambda_p(\chi) > 0 & \Longleftrightarrow \text{(1)} \\
\lambda_p(\chi) > 1 & \Longleftrightarrow \text{(1) and  (2)} \\
\lambda_p(\chi) > 2 & \Longleftrightarrow \text{(1), (2) and (3).} 
\end{align*}
\label{lambda}
\end{cor}
The above results can be used to compute $\lambda$-invariants of Dirichlet characters. We propose two alternative sets of parameters:\\

\noindent {\bf 1.}  For $k=1,2,\dots, p-1$, let $t_0 = (1+p)^{1-k} - 1$ and $b=(1+p)^{(1-p)+(1-k)} - 1$, e.g., choosing $k=1$ gives $t_0 = 0$ and $b=(1+p)^{1-p}-1$.  We have $v_p(b-t_0)=v_p((1+p)^{1-p}-1)=1$ and the special values are
\begin{align}
\mathcal{F}_{\chi}(t_0) & =   L_p(1-k,\chi)= - (1-\chi\omega^{-k}(p)p^{k-1}) \frac{B_{k,\chi \omega^{-k}}}{k} , \label{F0} \\
\mathcal{F}_{\chi}'(t_0)& =  \frac{ L_p'(1-k,\chi)} { (1+p)^{1-k} \log_p ( 1+p) } , \label{FD0} \\
\mathcal{F}_{\chi}( b ) &= L_p ( 2-p-k, \chi) = - (1-\chi\omega^{-k}(p) p^{p+k-2}) \frac{B_{p+k-1,\chi \omega^{-k}}}{p+k-1} \ .  \label{FDD0} 
\end{align}

\noindent {\bf 2.} For $k=1,2,\dots, p-1$, let $t_0 = p$ and $b=(1+p)^{1-k} - 1$, e.g., choosing $k=1$ yields $b=0$. We have $v_p(b-t_0)=v_p(-kp)=1$ and
\begin{align}
\mathcal{F}_{\chi}(p) & =   L_p(1,\chi) ,  \label{G0} \\
\mathcal{F}_{\chi}'(p) & =  \frac{ L_p'(1,\chi)} { (1+p) \log_p ( 1+p) } .
 \label{GD0} 
\end{align}

For a given character $\chi = \theta\omega^{i}$, we would usually choose $k=i$ for $i=1,\dots,p-2$ and $k=p-1$ for $i=0$. Then $\chi \omega^{-k}=\theta$, which simplifies the computations of the generalized Bernoulli numbers.

We also need explicit formulas for $L_p(1,\chi)$, $L_p'(1,\chi)$, $L_p'(1-k,\chi)$. For the  derivative of $L_p(s,\chi)$ at $s=0$ we recall a result of Ferrero-Greenberg:

\begin{prop} (\cite{FeGr} Prop. 1) Let $d$ be the conductor of $\chi \omega^{-1}$. If $p \nmid d$ then
\begin{align*} L_p'(0,\chi) = \left ( \sum_{a=1}^d \chi \omega^{-1} (a) \log_p \Gamma_p \left(\frac{a}{d}\right)  \right) + ( 1 - \chi \omega^{-1} (p) ) B_{1,\chi\omega^{-1}} \log_p (d) ,
\end{align*}
where $\Gamma_p$ is the $p$-adic gamma function. 
\label{FerrGreen}
\end{prop}

If  $\chi \omega^{-1}(p)=1$ then one has a simpler formula (see \cite{GrDa}):
\begin{align}   L_p'(0,\chi) = \sum_{a=1}^d \chi \omega^{-1} (a) \log_p \Gamma_p \left(\frac{a}{d}\right)  = 
  \sum_{a \in (\ZZ /d \ZZ)^*/<p> } \chi \omega^{-1} (a) \log_p \prod_{n=0}^{F-1}  \Gamma_p \left( \left\{ \frac{ap^n}{d} \right\} \right), 
  \label{FeGrForm}
  \end{align}
where $F$ is the multiplicative order of $p \mod d$ and $\{ x \}$ denotes the fractional part of $x$.
The latter product of $p$-adic gamma values  is related to Gauss and Jacobi sums by a theorem of Gross and Koblitz \cite{GrKo}, as explained in the following.

Let $q=p^F$ and let $\phi : \Bf_q^* \rightarrow \CC^{\times}$ be a character  of order $d$. Let $\psi: \Bf_q \rightarrow \CC^{\times}$ be the homomorphism
$\psi(x) = \zeta_p^{\text{Tr}(x)}$ where $\zeta_p = e^{2\pi i / p}$ and where $\text{Tr}: \Bf_q \rightarrow \Bf_p$ is the trace map. Then the Gauss sum $g(\phi)$ of $\phi$ is defined by
$$ g(\phi) = - \sum_{x \in \Bf_q^*} \phi(x) \psi(x). $$

 Now Theorem 1.7 of \cite{GrKo} gives for $a \in \{1, \dots, d-1 \}$:
 $$ \log_p g(\phi^{-a}) =  \log_p \prod_{n=0}^{F-1}  \Gamma_p \left( \left\{ \frac{ap^n}{d} \right\} \right) $$

 Let $J(\phi_1, \phi_2) = - \sum_{a_1+a_2=1} \phi_1(a_1) \phi_2(a_2)$ be 
 the Jacobi sum of the characters $\phi_1, \phi_2$ over $\Bf_q$.  If $\phi_1 \phi_2 \neq {\bf 1}$ then $J(\phi_1, \phi_2) = \frac{g(\phi_1) g(\phi_2)}{g(\phi_1 \phi_2)}$.
Furthermore, $g(\phi_1) g(\phi_1^{-1}) = \phi_1(-1) q$. Now suppose that $\gcd(a,d)=1$. We can express $g(\phi^{-a})^2$,  $g(\phi^{-a})^3$, $\dots$, $g(\phi^{-a})^d$  in terms of Jacobi sums and obtain
  $$  g(\phi^{-a})^d = \phi^{-a}(-1) q J(\phi^{-a}, \phi^{-a}) J(\phi^{-a}, \phi^{-2a}) \cdots J(\phi^{-a}, \phi^{-(d-2)a}) =: J(a, \phi) .$$
The above formula (\ref{FeGrForm}) for $L_p'(0,\chi)$ (in the case $\chi\omega^{-1}(p)=1$) then becomes

 \begin{align}
  L_p'(0,\chi)  & = \sum_{a \in (\ZZ /d \ZZ)^*/<p> } \chi\omega^{-1}(a) \log_p g(\phi^{-a}) =  \frac{1}{d} \sum_{a \in (\ZZ /d \ZZ)^*/<p> } \chi\omega^{-1}(a)  \log_p J(a, \phi) .
 \label{GrKoForm}
  \end{align}
  
We shall also need a formula when $p$ divides the conductor of $\chi \omega^{-1}$.
For this we refer to a result of Washington (\cite{Wa0} formulas $(2_p)$ and $(4_p)$, see also \cite{lang} Ch. 17, Theorem 4.2):
\begin{prop} Let $F$ be any multiple of $p$ and the conductor of $\chi$. Then:
$$ L_p'(0,\chi) = \left( \sum_{\substack{a=1 \\ p \nmid a}}^F \chi \omega^{-1} (a) G_p \left( \frac{a}{F} \right)  \right) + ( 1 - \chi \omega^{-1} (p) ) B_{1,\chi\omega^{-1}} \log_p (F)$$
where $G_p$ is Diamond's $p$-adic log gamma function defined by
$$G_p (X) = \left(X- \frac{1}{2} \right) \log_p X - X + \sum_{j=2}^{\infty} \frac{B_j}{j(j-1)} X^{1-j} $$
for $|X|_p > 1$ (see \cite{Dia}).
\label{formwash}
\end{prop}
There is also an explicit formula for $L_p'(s,\chi)$ at negative integers $s=1-k$ 
involving the values of $\chi \omega^{-k}$. We refer to \cite{Wa0} for details.
For $s=1$ we obtain the following formulas:

\begin{prop} Let $\chi  \neq {\bf 1}$ be a non-trivial Dirichlet character, and let $F$ be any multiple of $p$ and the conductor of $\chi$. Then
\begin{align*}
 L_p(1,\chi) & =  \sum_{\substack{a=1 \\ p \nmid a}}^F \chi(a) \left( - \frac{1}{F} \log_p(a) + \frac{1}{2a} + \sum_{j=2}^{\infty} B_j \frac{F^{j-1}}{j a^j} \right)  , \\
 L_p'(1,\chi) & =  \sum_{\substack{a=1 \\ p \nmid a}}^F \chi(a) \left( \frac{\log_p(a)^2}{2F}   -  \frac{\log_p(a)}{2a} - \sum_{j=2}^{\infty} B_j \frac{F^{j-1}}{j a^j} \left(\log_p(a) - 
 H_{j-1} \right) \right)  ,
 \end{align*}
 \label{format1}
 where $H_{j-1} = \displaystyle\sum_{k=1}^{j-1} \frac{1}{k}$ denotes the $(j-1)$-th harmonic number.
 \end{prop}
\begin{proof} We begin with Washington's formula for the $p$-adic $L$-function (see \cite{Wa2} Theorem 5.11):
\begin{align} L_p(s,\chi) = \frac{1}{F} \frac{1}{s-1} \sum _{\substack{a=1 \\ p \nmid a}}^F \chi(a) \langle a \rangle^{1-s} \sum_{j=0}^{\infty} \binom{1-s}{j} B_j \frac{F^j}{a^j} 
\label{washform}
\end{align}
To compute  the coefficients $a_0 = L_p(1,\chi)$ and $a_1 = L_p'(1,\chi)$ in the expansion $$L_p(s,\chi) = a_0 + a_1 (s-1) + a_2 (s-1)^2 + \dots ,$$ we adapt a  
method in  the proof of \cite{Wa2} Theorem 5.12, refining it to obtain exact values rather than working modulo $p$. See also Exercise 5.4 from ibid. for the formula on $L_p(1,\chi)$.

For integers $a$ coprime to $p$, we have
\begin{align} \langle a \rangle^{1-s} = 1+ (1-s) \log_p(a) + \frac{(1-s)^2}{2} (\log_p(a))^2 + \dots 
\label{form1}
\end{align}
Higher-order terms do not contribute to $a_0$ or $a_1$. Similarly, expanding the binomial sum gives
\begin{align}
 \sum_{j=0}^{\infty} \binom{1-s}{j} B_j \frac{F^j}{a^j} = 1 + (1-s) ( - \frac{1}{2} ) \frac{F}{a} + \sum_{j=2}^{\infty} \frac{(1-s) (-s) \cdots (2-j-s) }{j!} B_j \frac{F^j}{a^j}  .
 \label{form2}
 \end{align}
 The expansion for the sum on the right-hand side at $s=1$ is
$$ (1-s) \left( \sum_{j=2}^{\infty} (-1) \frac{1}{j} B_j \frac{F^j}{a^j} \right) + (1-s)^2 \left(  \sum_{j=2}^{\infty}  \frac{H_{j-1}}{j}  B_j \frac{F^j}{a^j} \right) + \dots .$$
To derive this, let $f(s) =  (1-s)(-s) (-1-s) \cdots (2-j-s)$. For even $j \geq 2$, we compute the expansion of $f(s)$ at $s=1$. The coefficient of $(1-s)$ is 
$(-1) \cdots (1-j) = (-1) (j-1)!$. The coefficient of $(1-s)^2$ is
\begin{align*}
 (-2)(-3) \cdots (1-j) +  (-1)(-3) \cdots (1-j) + \dots +  (-1)(-2) \cdots (2-j) \\ = \sum_{k=1}^{j-1} \frac{(-1)(-2) \cdots (1-j)}{-k} = H_{j-1} \cdot (j-1)! \ ,
 \end{align*}
 where $H_{j-1}$ is the $(j-1)$-th harmonic number.
 Multiplying the expansions (\ref{form1}) and (\ref{form2}) yields
 \begin{align*} 1+ & (s-1)  \left(- \log_p(a) + \frac{F}{2a} + \sum_{j=2}^{\infty} B_j \frac{F^{j}}{j a^j} \right) \\
 + & (s-1)^2   \left( \frac{\log_p(a)^2}{2}   -  \frac{F \log_p(a)}{2a} - \log_p(a) \sum_{j=2}^{\infty} B_j \frac{F^{j}}{j a^j} +  
\sum_{j=2}^{\infty} \frac{H_{j-1}}{j} B_j \frac{F^{j}}{a^j} \right) + \cdots
 \end{align*}
 Substituting this into (\ref{washform}) gives the claimed formulas.
\end{proof}

In the expansion $L_p(s,\chi) = a_0 + a_1 (s-1) + a_2 (s-1)^2 + \dots$, the coefficients $a_i$ satisfy
$p \mid a_i$ for all $i\geq 1$ (see \cite{Wa2} Theorem 5.11). This yields the congruences
\begin{align}
L_p(1-p, \chi) & = a_0 + a_1 (-p) + a_2 (-p)^2 + \dots \equiv a_0 = L_p(1,\chi) \mod p^2 ,
 \label{FDDD0} \\
 L_p'(1-p, \chi) & = a_1 + 2 a_2 (-p) + 3 a_3 (-p)^2 + \dots \equiv a_1 = L_p'(1, \chi) \mod p^2 .
\label{GDDD0} 
\end{align}

\section{The rank one case}\label{Sect3}

The interpolation factor in (\ref{InterpFormula}) -- the Euler factor  of $L(s, \chi \omega^{-1})$ at $p$ --  can cause $L_p(s,\chi)$ to have a zero at $s=0$. Ferrero and Greenberg (see  \cite{FeGr}) established  that $L_p(s,\chi)$ can have at most a simple zero at $s=0$.

\begin{prop}
If $\chi\omega^{-1}(p)=1$, then $L_p(s,\chi)$  has a trivial zero at $s=0$, in which case $\mathcal{F}_{\chi}(T)$ also has a trivial zero at $T=0$. This happens
precisely  if $\chi = \theta \omega$, where $\theta$ is an odd character with conductor prime to $p$, and $\theta(p) = \chi \omega^{-1} (p) =1$.
In this situation, we have $\lambda_p(\chi) \geq 1$.
\label{condzero}
\end{prop}
\begin{proof} This follows from the interpolation property (\ref{InterpFormula}) for $n=1$. We have $ \chi \omega^{-1} (p) =\theta(p)=1$, which implies that the conductor of $\theta$ is prime to $p$. Note that $\theta$ is odd since we assumed that $\chi$ is even.
\end{proof}

\begin{exm*} Let $\theta$ be the non-trivial character associated with a quadratic extension $K/\QQ$. Then $L_p(s,\theta \omega)$ has a trivial zero at $s=0$ if and only if $\theta$ is odd and $\theta(p)=1$. This is equivalent to $K$ being imaginary quadratic and $p$ being a prime which splits in $K$.
\end{exm*}

We now give explicit criteria to distinguish between the cases $\lambda_p(\chi) = 1$, $\lambda_p(\chi) =2$ and $\lambda_p(\chi) \geq 3$.

\begin{thm} Let $\chi$ be a non-trivial even Dirichlet character of the first kind and suppose that $\chi \omega^{-1} (p) = 1 $, so that $L_p(s,\chi)$ has a trivial zero at $s=0$.
\begin{enumerate}
\item[(i)] We have $\lambda_p(\chi) > 1$ if and only if one of the following equivalent conditions is satisfied:
\begin{enumerate}
\item[(a)] $v_p ( L_p'(0,\chi) ) > 1$
\item[(b)] $v_p(L_p(1,\chi)) > 1 $ 
\item[(c)] $ v_p ( B_{p, \chi \omega^{-1}}) > 2  $
\end{enumerate}
\item[(ii)] Assume that $\lambda_p(\chi) > 1$. 
Then $\lambda_p(\chi) > 2$ if and only if 
$$ v_p \left( - ( 1-  p^{p-1} )  \frac{B_{p, \chi \omega^{-1}}} {p}  - \frac{ L_p'(0,\chi)}{\log_p(1+p)} ( (1+p)^{1-p} -1) \right) > 2 .$$
An equivalent condition is that
$$ v_p \left( \frac{ L_p'(1,\chi)} { (1+p) \log_p ( 1+p) }\, p - L_p(1,\chi) \right) > 2 . $$

\end{enumerate}
\label{theorem-trivial}
\end{thm}

\begin{proof} 
Let $b=(1+p)^{1-p}-1$ and use the first set of parameters given in Section \ref{Sect2}. Then $v_p(b)=1$, and by (\ref{FDD0}), 
$\mathcal{F}_{\chi}(b) = ( 1- \chi \omega^{-1}(p) p^{p-1} )  \frac{B_{p, \chi \omega^{-1}}} {p} =  ( 1- p^{p-1} )  \frac{B_{p, \chi \omega^{-1}}} {p}$. 
From (\ref{FDDD0}) we obtain $\mathcal{F}_{\chi}(b) \equiv L_p(1,\chi) \mod p^2$.
Now part (i) of the assertion follows from Corollary  \ref{lambda} and equation (\ref{FD0}). Part (ii) follows from Corollary \ref{lambda} and equations (\ref{FD0}), (\ref{G0}), (\ref{GD0}).
For the equivalent condition, we use the second set of parameters.
\end{proof}

\begin{rem*} All of the above special values can be efficiently computed.  The value $L_p'(0,\chi)$ can be computed using either formula (\ref{FeGrForm}) or (\ref{GrKoForm}).
 Formulas for $L_p(1,\chi)$ and $L_p'(1,\chi)$ are given in Proposition \ref{format1}.
 
\end{rem*}

\begin{exm*} Let $p \equiv 1 \mod 4$ and let $\chi = \chi_{-4} \omega$, where $\chi_{-4}$ denotes the odd Dirichlet character of conductor $4$ associated with the imaginary quadratic field $\QQ(i)$. Then $\chi_{-4}(p)=1$ and
the $p$-adic $L$-function $L_p(s,\chi)$ has a trivial zero at $s=0$.
Let $E_n$ be the $n$-th {\em Euler number}  defined by
$$ \sum_{n=0}^{\infty} E_n \frac{x^n}{n!} = \frac{2}{e^x + e^{-x}} .$$
We observe that the generating function for $B_{n,\chi_{-4}}$ is given by
$$ \frac{ xe^x - xe^{3x}}{e^{4x}-1} = x \frac{e^x( 1-e^{2x})}{(e^{2x}+1)(e^{2x}-1)} = - x \frac{1}{e^x + e^{-x}} .$$
Hence $\frac{E_{n-1}}{2} = - \frac{B_{n,\chi_{-4}}}{n}$ for all $n \geq 1$.
By Theorem \ref{theorem-trivial} (i), $\lambda_p ( \chi) = \lambda_p (\QQ(i)) > 1$ if and only if $v_p (B_{p,\chi_{-4}}) >2  \Longleftrightarrow v_p(E_{p-1}) > 1$. 

Similarly, consider the odd Dirichlet character $\chi_{-3}$ of order $3$. Now let $p \equiv 1 \mod 3$ and $\chi = \chi_{-3} \omega$. Here we consider the $n$-th {\em Glaisher I-number} $G_n$ defined by
$$ \sum_{n=0}^{\infty} G_n \frac{x^n}{n!} = \frac{3/2}{1+e^x + e^{-x}} .$$
The generating function for $B_{n,\chi_{-3}}$ is given by
$$ \frac{ xe^x - xe^{2x}}{e^{3x}-1} = x \frac{e^x( 1-e^{x})}{(e^{x}-1)(e^{2x}+e^x + 1)} = - x \frac{1}{1+e^x + e^{-x}} .$$
We obtain $\frac{2}{3} G_{n-1}= - \frac{B_{n,\chi_{-3}}}{n}$ and hence $\lambda_p ( \chi) = \lambda_p (\QQ(\sqrt{-3})) > 1$ if and only if $v_p (B_{p,\chi_{-3}}) >2  \Longleftrightarrow v_p(G_{p-1}) > 1$.  These results on Euler and Glaisher numbers were also established by other methods in \cite{Sto}.
\end{exm*}

\begin{cor} Let $\chi = \theta \omega$ where $\theta$ is odd of conductor prime to $p$, and assume that $\theta(p)=1$. Let $K$ be the totally imaginary field of degree $2r_2$ cut out by the character $\theta$. Then the $p$-adic $L$-functions $L_p( \theta^{i} \omega ,s)$, where $ i =1, 3, \dots, 2r_2 -1$, have a trivial zero at $s=0$.
\end{cor}
\begin{proof}  
This follows from $\theta^i (p) = 1$ for all $i \in \NN$ and Proposition \ref{condzero}.
\end{proof}

The assumption $\theta(p)=1$ implies that $p$ splits completely in $K$.
Note that $\theta, \theta^3, \dots, \theta^{2r_2-1}$ are the odd characters associated with $K$. 
Now let $\lambda_p^- (K) = \sum_{i=1}^{r_2} \lambda_p(\theta^{2i-1} \omega)$. Then the above
Corollary yields $\lambda_p^-(K) \geq r_2$. 

The following Proposition gives the corresponding fact about the growth of the $p$-primary part of the ideal class groups in the cyclotomic extension of $K$. In fact, this is known for more general number fields and $\Zp$-extensions: 

\begin{prop} Let $K$ be a finite extension of $\QQ$ and assume that $p$ splits completely in $K/\QQ$. Let $r_2$ be the number of pairs of complex embeddings of $K$.
Assume that $K_{\infty}$ is a $\Zp$-extension of $K$ such that every prime above $p$ is ramified in $K_{\infty}/K$, e.g., that $K_{\infty}$ is the cyclotomic extension of $K$. Then the Iwasawa $\lambda_p$-invariant of $K$ (defined by the growth of the $p$-primary ideal class groups in the $\Zp$-extension $K_{\infty}/K$) satisfies
$$\lambda_p(K) \geq r_2 .$$ 
\end{prop}
\begin{proof} 
This can be shown using Galois theory and class field theory (see \cite{Gr} Prop. 2.2).
\end{proof}

\section{The rank zero case}\label{Sect4}

We now address the scenario where $L_p(0,\chi)$ does not vanish. The following theorem provides criteria for the cases $\lambda_p(\chi)=0$, $\lambda_p(\chi) = 1$, $\lambda_p(\chi) =2$ and $\lambda_p(\chi) \geq 3$. Notably, the result in part (i) is well known.

\begin{thm} Let $\chi$ be a non-trivial even character of the first kind and suppose that $\chi \omega^{-1} (p) \neq 1 $. Then:
\begin{enumerate}
\item[(i)] We have $\lambda_p(\chi) > 0$ if and only if $v_p(B_{1,\chi \omega^{-1}}) = v_p (L_p(0,\chi) ) > 0$.
\item[(ii)] Assume that $\lambda_p(\chi) > 0$. Then $\lambda_p(\chi) > 1$ if one of the following equivalent conditions is satisfied:
\begin{enumerate}
\item[(a)]  $v_p (L_p'(0,\chi)) > 1$
\item[(b)] $ v_p (- (1-\chi\omega^{-1}(p)) B_{1,\chi \omega^{-1}} - L_p(1,\chi)  ) > 1  $
\item[(c)] $ v_p \left( -( 1- \chi \omega^{-1}(p) p^{p-1} )  \frac{B_{p, \chi \omega^{-1}}} {p} + (1-\chi\omega^{-1}(p)) B_{1,\chi \omega^{-1}} \right) > 1  $
\end{enumerate}

\item[(iii)] Assume that $\lambda_p(\chi) > 1$. 
Then $\lambda_p(\chi) > 2$ if and only if 
$$
\begin{multlined}
 v_p \bigg( - \left( 1- \chi \omega^{-1}(p) p^{p-1} \right)  \frac{B_{p, \chi \omega^{-1}}} {p}  + \left(1- \chi \omega^{-1}(p) \right) B_{1, \chi \omega^{-1}} \\
 -  \frac{ L_p'(0,\chi)}{\log_p(1+p)} \left( (1+p)^{1-p} -1 \right) \bigg) > 2.
 \end{multlined}
 $$
 An equivalent condition is that
$$ v_p \left( - (1-\chi\omega^{-1}(p)) B_{1,\chi \omega^{-1}} + \frac{ L_p'(1,\chi)} { (1+p) \log_p ( 1+p) }\, p - L_p(1,\chi) \right) > 2 . $$

\end{enumerate}
\label{theorem-general}
\end{thm}

\begin{proof}
As in the proof of Theorem \ref{theorem-trivial}, the assertions follow from Corollary \ref{lambda} and equations (\ref{F0}), (\ref{FD0}), (\ref{FDD0}), (\ref{G0}), (\ref{GD0}). 
\end{proof}

\begin{rem*} All of the above special values can be efficiently computed using Propositions \ref{FerrGreen}, \ref{formwash} and \ref{format1}.
\end{rem*}

\begin{rem*} Let us fix an odd character $\theta$, and let $K$ be the associated totally imaginary field. We know from Section \ref{Sect3} that $\lambda_p(\theta \omega) \geq 1$ for odd primes $p$ that split completely in $K/\QQ$, i.e., primes $p$ with $\theta(p)=1$. The number of such primes is infinite. One can use the $p$-adic valuation of  $B_{p,\theta}$ or the valuation of $L_p(1,\chi)$ to determine whether $\lambda_p(\theta \omega) =1$ or $\lambda_p(\theta \omega) > 1$ (see Theorem \ref{theorem-trivial} (i)).

On the other hand, consider the non-split unramified primes $p$, i.e., primes with $\theta(p) \neq 0, 1$. There are infinitely many such primes, but only finitely many have the property that $\lambda_p(\theta \omega) > 0$  since $v_p(B_{1,\theta})=0$ for almost all primes $p$.

\end{rem*}

\section{Twists by characters of $p$-power order } 
\label{Sect5}  
In this section, we explore how twists by characters of $p$-power order can be leveraged to compute $\lambda$-invariants.
 Dummit et al.\ \cite{Du} used characters $\psi$ of order $p^n$ to set up a system of equations and compute the coefficients of $\mathcal{F}_{\chi}(T)$  modulo $p^n$ (or modulo $p^{n-1}$, as applicable).
 We follow the approach by Ernvall and Metsänkylä  in \cite{ErMe} and show that the $p$-adic valuation of a {\em single value} $L_p(1-k,\chi \psi)$ suffices to compute $\lambda_p(\chi)$ if $n$ is large enough ($n=1$ or $n=2$ are usually sufficient).

\begin{prop} Let $K/\QQ_p$ be a finite extension with ring of integers $\mathcal{O}_K$ and ramification index $e$ at $p$. Let $n$ be a positive integer  and let 
$\pi_n $ be a uniformizer of $\QQ_p(\zeta_{p^n})$, 
e.g.,  $\pi_n = \zeta_{p^n} - 1$. 
Let $t_0 \in \mathcal{O}_K$ with $v_p(t_0) \geq 1$, e.g., $t_0=0$. 
Suppose that
$\mathcal{F}(T)  \in \mathcal{O}_K[[T]]$ is a power series satisfying $\mu_p(\mathcal{F})=0$. 
\begin{enumerate}
\item[(i)] We have $\lambda_p(\mathcal{F}) < \frac{(p-1)p^{n-1}}{e}$ if and only if  $v_p(\mathcal{F}(\pi_n )) < \frac{1}{e}$, and in this case
$$ \lambda_p(\mathcal{F}) = v_p(\mathcal{F}(\pi_n ))\cdot (p-1)p^{n-1} .\\ $$ 
\item[(ii)] Suppose that $\lambda_p(\mathcal{F}) > 0$. Then the equivalence 
$$\lambda_p(\mathcal{F}) < \frac{(p-1)p^{n-1}}{e} + 1 \Longleftrightarrow v_p(\mathcal{F}(\pi_n ) - \mathcal{F}(t_0)) <  \frac{1}{e} + \frac{1}{(p-1)p^{n-1}}$$
holds. In this case,
$$  \lambda_p(\mathcal{F})  = v_p(\mathcal{F}(\pi_n ) - \mathcal{F}(t_0)) \cdot (p-1)p^{n-1}   .$$
\item[(iii)] Suppose that $\lambda_p(\mathcal{F}) > 1$. Then  the equivalence 
$$\lambda_p(\mathcal{F}) < \frac{(p-1)p^{n-1}}{e} + 2 \Longleftrightarrow v_p(\mathcal{F}(\pi_n ) - \mathcal{F}(t_0) -   \mathcal{F}'(t_0)(\pi_n - t_0)) <  \frac{1}{e} + \frac{2}{(p-1)p^{n-1}}$$
holds. In this case,
$$  \lambda_p(\mathcal{F})  = v_p(\mathcal{F}(\pi_n ) - \mathcal{F}(t_0) -   \mathcal{F}'(t_0 ) (\pi_n - t_0) ) \cdot (p-1)p^{n-1}  .$$

\end{enumerate}

\label{secondkind}
\end{prop}

\begin{proof} The statement is true for $\lambda= \lambda_p(\mathcal{F})=0$, and we can assume that $\lambda   \geq 1$. 
Let $\mathcal{F}(T)=a_0 + a_1 (T-t_0) + \dots + a_{\lambda-1}(T-t_0)^{\lambda -1} + a_{\lambda} (T-t_0)^{\lambda} + \dots $. 
We extend the valuation $v_p$ to all finite extensions of $K$ inside a fixed algebraic closure. 
Since $v_p(t_0) > 0$,  
the characterization of $\lambda$ via valuations of the coefficients $a_j$ does not depend on the chosen center $t_0$.
Furthermore, $v_p(a_0), \dots , v_p(a_{\lambda-1}) \geq  \frac{1}{e}$ and $v_p(a_{\lambda})=0$. By assumption, $v_p(\pi_n) = v_p(\pi_n - t_0) = \frac{1}{(p-1)p^{n-1}}$. If $\lambda <  \frac{(p-1)p^{n-1}}{e}$ we have $v_p(\mathcal{F}(\pi_n)) = \frac{\lambda}{(p-1)p^{n-1}} < \frac{1}{e}$. Conversely, suppose that $\lambda \geq \frac{(p-1)p^{n-1}}{e}$.
Then  $v_p(\mathcal{F}(\pi_n )) \geq \frac{1}{e}$. 
This proves (i), and parts (ii) and (iii) can be shown similarly by subtracting the terms of degree $0$ and $1$  from $\mathcal{F}(T)$.
\end{proof}

\begin{rem*} \begin{enumerate}
\item We can increase $n$ in order to satisfy the condition on $\lambda_p(\mathcal{F})$, and thus compute the exact $\lambda$-invariant. In most cases, $n=1$ or $n=2$ are sufficient since large values of $\lambda_p(\chi)$ are very scarce. Obviously, one would choose $n$ as small as possible for explicit computations.

\item Let $\psi$ be a character of order $p^n$ where $n\geq 1$, and let $\zeta_{p^n}=\zeta_{\psi}=\psi(1+p)^{-1}$.
The obvious choices for $t_0$ and $\pi_n$ are as follows: for $k=1,\dots,p-1$, set $t_0=(1+p)^{1-k}-1$ and $\pi_n = \zeta_{p^n}(1+p)^{1-k}-1$. Then 
 $\mathcal{F}_{\chi}(t_0) = L_p(1-k,\chi)$ and  
 $\mathcal{F}_{\chi}(\pi_n) =  
L_p(1-k,\chi\psi)=- \frac{1}{k} B_{k,\chi\psi\omega^{-k}}$ (see formula (\ref{powerseries})).

\item Parts (ii) and (iii) of the above Proposition turn out to be useful for $\lambda >0$, small primes $p$, and $n=1$ or $n=2$. In our applications, the constant $\mathcal{F}(t_0)$ and the derivative  $\mathcal{F}'(t_0)$ can be efficiently computed.
\end{enumerate}
\end{rem*}

The Theorem below is a slight improvement of \cite{ErMe} through incorporating the values $L_p(1-k,\chi)$ and $L_p'(1-k,\chi)$. The assertions follow from Proposition \ref{secondkind}, where $t_0$ and $\pi_n$ are chosen as in above Remark (2.), and formulas (\ref{InterpFormula}), (\ref{F0}), (\ref{FD0}).

\begin{thm} Let $\chi = \theta \omega^i$ be a non-trivial Dirichlet character of the first kind, and let $e$ be the ramification index of $\mathcal{O}_{\chi}$. Define $k=i$ if $i=1,\dots,p-2$, and $k=p-1$ if $i=0$. Let $\psi$ be a character of the second kind of order $p^n$, where $n\geq 1$. 
\begin{enumerate}
\item[(i)] (See \cite{ErMe}). We have
$ \lambda_p(\chi) < \frac{(p-1)p^{n-1}}{e}   \Longleftrightarrow v_p(B_{k,\theta \psi} ) <  \frac{1}{e} $,
and if this  is satisfied then
$$  \lambda_p(\chi) = v_p(B_{k,\theta \psi} ) \cdot (p-1)p^{n-1} . $$

\item[(ii)]  Suppose that $\lambda_p(\chi)>0$. Then the equivalence 
$$\lambda_p(\chi) < \frac{(p-1)p^{n-1}}{e} + 1 \Longleftrightarrow v_p\left( - \frac{B_{k,\theta \psi}}{k} + (1-\theta(p)p^{k-1}) \frac{B_{k,\theta}}{k}\right) < \frac{1}{e} + \frac{1}{(p-1)p^{n-1}}$$
holds, and in this case, 
$$  \lambda_p(\chi) = v_p\left(- \frac{B_{k,\theta \psi}}{k} + (1-\theta(p)p^{k-1}) \frac{B_{k,\theta }}{k} \right) \cdot (p-1)p^{n-1} . $$

\item[(iii)] Suppose that $\lambda_p(\chi)>1$.  Then 
$$
\begin{multlined} \lambda_p(\chi) < \frac{(p-1)p^{n-1} }{e} + 2 \Longleftrightarrow \\  
v_p \left( - \frac{B_{k,\theta \psi}}{k} + (1-\theta(p)p^{k-1}) \frac{B_{k,\theta }}{k} - \frac{L_p'(1-k,\theta )}{\log_p(1+p)}\: (\zeta_{p^n}-1) \right) 
< \frac{1}{e} + \frac{2}{(p-1)p^{n-1}}
\end{multlined}
$$
holds, and in this case, $\lambda_p(\chi)$ is given by the valuation multiplied by $(p-1)p^{n-1}$.

\end{enumerate}
\label{thm-val}
\end{thm} 

Alternatively, we can also utilize the values of the $p$-adic $L$-functions associated with $\chi$ and $\chi\psi$ at $s=1$. The value and the derivative at $s=1$  can be efficiently computed using Proposition \ref{format1}.
The Theorem below follows from Proposition \ref{secondkind}, where we choose $t_0=p$ and $\pi_n = \zeta_{p^n}(1+p)-1$, and formulas (\ref{G0}) and (\ref{GD0}). We have $\mathcal{F}_{\chi}(p)=L_p(1,\chi)$ and $\mathcal{F}_{\chi}(\pi_n)=L_p(1,\chi \psi)$.

\begin{thm} Let $\chi $ be a non-trivial Dirichlet character of the first kind, and let $e$ be the ramification index of $p$ in $\mathcal{O}_{\chi}$. Let $\psi$ be a character of the second kind of order $p^n$, where $n\geq 1$. 
\begin{enumerate}
\item[(i)] We have
$ \lambda_p(\chi) < \frac{(p-1)p^{n-1}}{e}   \Longleftrightarrow v_p(L_p(1,\chi \psi)) <  \frac{1}{e} $,
and if this  is satisfied then
$$  \lambda_p(\chi) = v_p(L_p(1,\chi \psi) ) \cdot (p-1)p^{n-1} . $$

\item[(ii)] Suppose that $\lambda_p(\chi)>0$. Then the equivalence
$$
\begin{multlined}
\lambda_p(\chi) < \frac{(p-1)p^{n-1}}{e} + 1 \Longleftrightarrow \\ 
v_p\left( L_p(1,\chi \psi) - L_p(1,\chi)   \right) < \frac{1}{e} + \frac{1}{(p-1)p^{n-1}}
\end{multlined}
$$
holds, and in this case,
$$  \lambda_p(\chi) = v_p\left( L_p(1,\chi \psi) - L_p(1,\chi)   \right) \cdot (p-1)p^{n-1} . $$

\item[(iii)] Suppose that $\lambda_p(\chi)>1$. Then the equivalence
$$
\begin{multlined} \lambda_p(\chi) < \frac{(p-1)p^{n-1}}{e} + 2 \Longleftrightarrow \\ 
v_p \left( L_p(1,\chi \psi) - L_p(1,\chi)   - \frac{L_p'(1,\chi )}{\log_p(1+p)}\: (\zeta_{p^n}-1) \right) < \frac{1}{e} + \frac{2}{(p-1)p^{n-1}}
\end{multlined}
$$
holds, and in this case, $\lambda_p(\chi)$ is given by the valuation multiplied by $(p-1)p^{n-1}$.

\end{enumerate}
\label{thm1}
\end{thm} 

\begin{rem*} The above two theorems (in particular parts (i) and (ii)) turn out to be very useful for practical computations of $\lambda_p(\chi)$ (see Section \ref{Sect6} below). The running time is (as expected) linear in the conductor of the characters involved, and computations are quite fast on standard PCs for conductors $<10^6$.

\end{rem*}

\section{Numerical data}\label{Sect6}
In this section, we extend our previous work
by computing $\lambda$-invariants via newly derived formulas. We also collect new data for the rank one case supporting our prediction regarding the distribution of $\lambda$-invariants. 

For most computations, we use the special values of the first set of parameters (see Section \ref{Sect2}).  As before, assume $\chi = \theta \omega^i$ with $p$ not dividing the conductor of $\theta$.  In the rank one case, the $p$-th generalized Bernoulli number of $\theta = \chi \omega^{-1} $ and the derivative $L'_p(0,\chi)$  (modulo $p^3$) are sufficient to distinguish between the cases $\lambda=1$, $\lambda=2$, and $\lambda >2$. In the rank zero case, the value $B_{1,\theta}$ is also needed. Furthermore, to compute $\lambda$-invariants $\geq 3$, we use the $k$-th generalized Bernoulli numbers (where $k=1,\dots,p-1$) of $\theta$ and of $\theta \psi$ (where $\psi$ has order $p^n$), as described in part (ii) of Theorem \ref{thm-val}, starting with $n=1$ and increasing $n$ if necessary.

 We have cross-checked the results by using different formulas, in particular with the second set of parameters, characters of the second kind, the $p$-adic expansions in \cite{KnWa}, and also calculations by other authors, e.g., \cite{Du}. The numerical computations were done with SageMath and the code is available from the website
\url{https://github.com/knospe/iwasawa}.

\subsection{Distribution of $\lambda$-invariants in the rank one case}
\label{Sect41}
First, we provide numerical data for the rank one case and fixed primes $p$. In this situation, the predicted distribution of $\lambda$-values (see \cite{DeKn}) should be shifted by $1$. 

Using the above formulas, we computed the $\lambda$-invariants of $\chi = \theta \omega$, where $\theta$ is odd of conductor $<100,000$ and $\theta(p)=1$. The first row gives the predicted distribution and the second row contains the number $N$ of tested characters as well as the computed proportions. For $p=3$ and $\ord(\theta)=2$ we still see a larger deviation from the prediction (compare \cite{DeKn}), but the difference becomes smaller as we compute $\lambda$-invariants of characters with larger conductor. \\

\noindent $p=3$, $\ord(\theta)=2$\\
\begin{tabular}{|l|l|l|l|l|l|l|l|} \hline
& $\lambda=1$ & $\lambda=2$ & $\lambda=3$ & $\lambda=4$ & $\lambda=5$ & $\lambda=6$ & $\lambda\geq 7$ \\ \hline
predicted & $0.5601$ & $0.2801$ & $0.1050$ & $0.0364$ & $0.0123$ & $0.0041$ & $0.0020$ \\ \hline
$N=11404$ & $0.6121$ & $0.2604$ & $0.0835$ & $0.0296$ & $0.0092$ & $0.0035$ & $0.0017$ \\ \hline
\end{tabular}\\

\noindent $p=3$, $\ord(\theta)=4$\\
\begin{tabular}{|l|l|l|l|l|l|l|l|} \hline
& $\lambda=1$ & $\lambda=2$ & $\lambda=3$ & $\lambda=4$ & $\lambda=5$ & $\lambda=6$ & $\lambda\geq 7$ \\ \hline
predicted & $0.8766$ & $0.1096$ & $0.0123$ & $0.0014$ & $0.0002$ & $0.0000$ & $0.0000$ \\ \hline
$N=17732$ & $0.8901$ & $0.0976$ & $0.0108$ & $0.0014$ & $0.0001$ & $0.0000$ & $0.0000$ \\  \hline
\end{tabular}\\

\noindent $p=5$, $\ord(\theta)=2$\\
\begin{tabular}{|l|l|l|l|l|l|l|l|} \hline
& $\lambda=1$ & $\lambda=2$ & $\lambda=3$ & $\lambda=4$ & $\lambda=5$ & $\lambda=6$ & $\lambda\geq 7$ \\ \hline
predicted & $0.7603$ & $0.1901$ & $0.0396$ & $0.0080$ & $0.0016$ & $0.0003$ & $0.0001$\\ \hline
$N=12667$ & $0.7782$ & $0.1768$ & $0.0367$ & $0.0072$ & $0.0007$ & $0.0003$ & $0.0001$  \\  \hline
\end{tabular}\\

\noindent $p=5$, $\ord(\theta)=4$\\
\begin{tabular}{|l|l|l|l|l|l|l|l|} \hline
& $\lambda=1$ & $\lambda=2$ & $\lambda=3$ & $\lambda=4$ & $\lambda=5$ & $\lambda=6$ & $\lambda\geq 7$ \\ \hline
predicted & $0.7603$ & $0.1901$ & $0.0396$ & $0.0080$ & $0.0016$ & $0.0003$ & $0.0001$\\ \hline
$N=12679$ & $0.7662$ & $0.1834$ & $0.0414$ & $0.0072$ & $0.0013$ & $0.0001$ & $0.0003$ \\ \hline
\end{tabular}\\

\noindent $p=5$, $\ord(\theta)=6$\\
\begin{tabular}{|l|l|l|l|l|l|l|l|} \hline
& $\lambda=1$ & $\lambda=2$ & $\lambda=3$ & $\lambda=4$ & $\lambda=5$ & $\lambda=6$ & $\lambda\geq 7$ \\ \hline
predicted & $0.9584$ &  $0.0399$ &  $0.0016$ &  $0.0001$ & $0.0000$ &  $0.0000$ & $0.0000$ \\  \hline
$N=37291$ & $0.9578$ & $0.0407$ & $0.0015$ & $0.0001$ & $0.0000$ & $0.0000$ & $0.0000$ \\ \hline
\end{tabular}\\

\noindent $p=7$, $\ord(\theta)=2$\\
\begin{tabular}{|l|l|l|l|l|l|l|l|} \hline
& $\lambda=1$ & $\lambda=2$ & $\lambda=3$ & $\lambda=4$ & $\lambda=5$ & $\lambda=6$ & $\lambda\geq 7$ \\ \hline
predicted & $0.8368$ & $0.1395$ &  $0.0203$ &  $0.0029$ &  $0.0004$ &  $0.0001$ & $0.0000$ \\ \hline
$N=13299$ & $0.8419$ & $0.1341$ & $0.0207$ & $0.0030$ & $0.0002$ & $0.0000$ & $0.0001$\\ \hline
\end{tabular}\\

\noindent $p=7$, $\ord(\theta)=4$\\
\begin{tabular}{|l|l|l|l|} \hline
& $\lambda=1$ & $\lambda=2$ & $\lambda \geq 3$ \\ \hline
predicted & $0.9792$ & $0.0204$ & $0.0004$ \\ \hline
 $N=21244$ & $0.9797$ & $0.0197$ & $0.0007$ \\ \hline
\end{tabular}\\

\noindent $p=7$, $\ord(\theta)=6$\\
\begin{tabular}{|l|l|l|l|l|l|l|l|} \hline
& $\lambda=1$ & $\lambda=2$ & $\lambda=3$ & $\lambda=4$ & $\lambda=5$ & $\lambda=6$ & $\lambda\geq 7$ \\ \hline
predicted & $0.8368$ & $0.1395$ &  $0.0203$ &  $0.0029$ &  $0.0004$ &  $0.0001$ & $0.0000$ \\ \hline
$N=27049$ & $0.8360$ & $0.1419$ & $0.0187$ & $0.0029$ & $0.0003$ & $0.00004$ & $0.00004$\\ \hline
\end{tabular}\\

\noindent $p=11$, $\ord(\theta)=2$\\
\begin{tabular}{|l|l|l|l|} \hline
 & $\lambda=1$ & $\lambda=2$ & $\lambda \geq 3$ \\ \hline
predicted & $0.9008$ & $0.0901$ & $0.0091$ \\ \hline
$N=13931$ & $0.9058$ & $0.0847$ & $0.0095$ \\ \hline

\end{tabular}\\

\noindent $p=11$, $\ord(\theta)=4$\\
\begin{tabular}{|l|l|l|l|} \hline
 & $\lambda=1$ & $\lambda=2$ & $\lambda \geq 3$ \\ \hline
predicted & $0.9917$ & $0.0083$ & $0.0001$ \\ \hline
$N=23252$ & $0.9895$ & $0.0104$ & $0.0001$\\ \hline
\end{tabular}\\

\noindent $p=11$, $\ord(\theta)=6$\\
\begin{tabular}{|l|l|l|l|} \hline
 & $\lambda=1$ & $\lambda=2$ & $\lambda \geq 3$ \\ \hline
predicted & $0.9917$ & $0.0083$ & $0.0001$ \\ \hline
$N=41478$ & $0.9920$ & $0.0079$ & $0.0002$\\ \hline
\end{tabular}\\

\noindent $p=13$, $\ord(\theta)=2$\\
\begin{tabular}{|l|l|l|l|} \hline
 & $\lambda=1$ & $\lambda=2$ & $\lambda \geq 3$ \\ \hline
predicted & $0.9172$ & $0.0764$ & $0.0064$ \\ \hline
$N=14110$ & $0.9199$ & $0.0728$ & $0.0073$\\  \hline

\end{tabular}\\

\noindent $p=13$, $\ord(\theta)=4$\\
\begin{tabular}{|l|l|l|l|} \hline
 & $\lambda=1$ & $\lambda=2$ & $\lambda \geq 3$ \\ \hline
predicted & $0.9172$ & $0.0764$ & $0.0064$  \\ \hline
$N=18806$ & $0.9210$ & $0.0726$ & $0.0064$\\ \hline
\end{tabular}\\

\noindent $p=13$, $\ord(\theta)=6$\\
\begin{tabular}{|l|l|l|l|} \hline
 & $\lambda=1$ & $\lambda=2$ & $\lambda \geq 3$ \\ \hline
predicted & $0.9172$ & $0.0764$ & $0.0064$  \\ \hline
$N=36440$ & $0.9189$ & $0.0743$ & $0.0068$ \\ \hline
\end{tabular}\\

\noindent $p=17$, $\ord(\theta)=2$\\
\begin{tabular}{|l|l|l|l|} \hline
 & $\lambda=1$ & $\lambda=2$ & $\lambda \geq 3$ \\ \hline
predicted & $0.9377$ & $0.0586$ & $0.0037$ \\ \hline
$N=14352$ & $0.9392$ & $0.0575$ & $0.0033$ \\ \hline
\end{tabular}\\

\noindent $p=19$, $\ord(\theta)=2$\\
\begin{tabular}{|l|l|l|l|} \hline
 & $\lambda=1$ & $\lambda=2$ & $\lambda \geq 3$ \\ \hline
predicted &  $0.9446$ & $0.0525$ & $0.0029$ \\ \hline
$N=14441$ & $0.9463$ & $0.0504$ & $0.0033$ \\ \hline
\end{tabular}\\

\noindent $p=23$, $\ord(\theta)=2$\\
\begin{tabular}{|l|l|l|l|} \hline
 & $\lambda=1$ & $\lambda=2$ & $\lambda \geq 3$ \\ \hline
predicted &  $0.9546$ & $0.0434$ & $0.0020$ \\ \hline
$N=14569$ &  $0.9552$ & $0.0428$ & $0.0020$\\ \hline
\end{tabular}\\

\noindent $p=29$, $\ord(\theta)=2$\\
\begin{tabular}{|l|l|l|l|} \hline
 & $\lambda=1$ & $\lambda=2$ & $\lambda \geq 3$ \\ \hline
predicted &  $0.9643$ & $0.0344$ & $0.0012$ \\ \hline
$N=14695$ & $0.9666$ & $0.0324$ & $0.0010$ \\ \hline
\end{tabular}\\

\noindent $p=31$, $\ord(\theta)=2$\\
\begin{tabular}{|l|l|l|l|} \hline
 & $\lambda=1$ & $\lambda=2$ & $\lambda \geq 3$ \\ \hline
predicted &  $0.9667$ & $0.0322$ & $0.0011$ \\ \hline
 $N=14728$ & $0.9659$ & $0.0328$ & $0.0013$\\ \hline
\end{tabular}\\

\noindent $p=31$, $\ord(\theta)=6$\\
\begin{tabular}{|l|l|l|l|} \hline
 & $\lambda=1$ & $\lambda=2$ & $\lambda \geq 3$ \\ \hline
predicted &  $0.9667$ & $0.0322$ & $0.0011$ \\ \hline
 $N=41770$ & $0.9672$ & $0.0318$ & $0.0011$\\ \hline
\end{tabular}\\

\noindent $p=37$, $\ord(\theta)=2$\\
\begin{tabular}{|l|l|l|l|} \hline
 & $\lambda=1$ & $\lambda=2$ & $\lambda \geq 3$ \\ \hline
predicted & $0.9722$ & $0.0270$ & $0.0008$   \\ \hline
$N=14789$ &  $0.9734$ & $0.0260$ & $0.0007$  \\ \hline
\end{tabular}\\

\noindent $p=41$, $\ord(\theta)=2$\\
\begin{tabular}{|l|l|l|l|} \hline
 & $\lambda=1$ & $\lambda=2$ & $\lambda \geq 3$ \\ \hline
predicted &  $0.9750$ & $0.0244$ & $0.0006$ \\ \hline
$N=14841$ & $0.9749$ & $0.0245$ & $0.0007$  \\ \hline
\end{tabular}\\

\noindent $p=43$, $\ord(\theta)=2$\\
\begin{tabular}{|l|l|l|l|} \hline
 & $\lambda=1$ & $\lambda=2$ & $\lambda \geq 3$ \\ \hline
predicted & $0.9762$ & $0.0232$ & $0.0006$ \\ \hline
$N=14856$ & $0.9768$ & $0.0228$ & $0.0004$ \\ \hline
\end{tabular}\\

\noindent $p=43$, $\ord(\theta)=6$\\
\begin{tabular}{|l|l|l|l|} \hline
 & $\lambda=1$ & $\lambda=2$ & $\lambda \geq 3$ \\ \hline
predicted & $0.9762$ & $0.0232$ & $0.0006$ \\ \hline
$N=46382$ & $0.9777$ & $0.0218$ & $0.0005$ \\ \hline
\end{tabular}\\

\noindent $p=47$, $\ord(\theta)=2$\\
\begin{tabular}{|l|l|l|l|} \hline
 & $\lambda=1$ & $\lambda=2$ & $\lambda \geq 3$ \\ \hline
predicted & $0.9783$ & $0.0213$ & $0.0005$ \\ \hline
$N=14884$ & $0.9789$ & $0.0204$ & $0.0007$ \\ \hline
\end{tabular}\\

\noindent $p=97$, $\ord(\theta)=2$\\
\begin{tabular}{|l|l|l|l|} \hline
 & $\lambda=1$ & $\lambda=2$ & $\lambda \geq 3$ \\ \hline
predicted & $0.9896$ & $0.0103$ & $0.0001$ \\ \hline
$N=15040$ & $0.9897$ & $0.0102$ & $0.0001$  \\ \hline
\end{tabular}

\subsection{Primes with $\lambda>1$ in the rank one case}
In this subsection, we identify small primes and fields of small degree for which the $\lambda$-invariant exceeds $1$ in the rank one case.
We fixed an odd character $\theta$ and utilized Theorem \ref{theorem-trivial} (i) to compute all small primes $p$ such that  $L_p(s,\theta \omega)$ has a trivial zero and $\lambda_p(\theta \omega)>1$. This was also done by Dummit et al.\ in \cite{Du} for imaginary quadratic characters. 
Below, we list the combinations of characters $\chi = \theta \omega$ of order $<10$, conductor $<100$, and primes $p<500$ having a trivial zero and satisfying $\lambda_p(\chi)>1$. For a fixed character, we expect (see \cite{Du, DeKn}) that the number of such primes $p$ with $p \leq X$ is $O(\log(\log(X))$. Furthermore, almost all of these primes should have residue class degree $f=1$, i.e.,  only a finite number is expected to have degree $f \geq 2$. \\

\noindent \begin{tabular}{|l|l|l|l|} \hline
order & conductor & field & primes  \\ \hline
$2$ & $3$ & $\QQ(\sqrt{-3})$ & $13 $, $181$ \\ \hline
$2$ & $11$ & $\QQ(\sqrt{-11})$ & $5 $  \\ \hline
$2$ & $19$ & $\QQ(\sqrt{-19})$ & $11 $  \\ \hline
$2$ & $24$ & $\QQ(\sqrt{-6})$ & $131 $  \\ \hline
$2$ & $31$ & $\QQ(\sqrt{-31})$ & $227 $  \\ \hline
$2$ & $35$ & $\QQ(\sqrt{-35})$ & $3$, $13$  \\ \hline
$2$ & $47$ & $\QQ(\sqrt{-47})$ & $3$, $17$, $157$  \\ \hline
$2$ & $51$ & $\QQ(\sqrt{-51})$ & $5 $  \\ \hline
$2$ & $52$ & $\QQ(\sqrt{-13})$ & $113 $  \\ \hline
$2$ & $56$ & $\QQ(\sqrt{-14})$ & $3 $  \\ \hline
$2$ & $71$ & $\QQ(\sqrt{-71})$ & $29 $  \\ \hline
$2$ & $83$ & $\QQ(\sqrt{-83})$ & $17 $, $41$  \\ \hline
$2$ & $84$ & $\QQ(\sqrt{-21})$ & $107 $, $173$  \\ \hline
$2$ & $88$ & $\QQ(\sqrt{-22})$ & $23$, $29$  \\ \hline

\end{tabular}

\bigskip

\noindent \begin{tabular}{|l|l|l|l|} \hline
order & conductor & field (polynomial) & prime\\ \hline
$4$ & $16$ & $x^4 + 4 x^2 + 2$ & $97 $ \\ \hline
$4$ & $29$ & $x^4 + x^3 + 4 x^2 + 20 x + 23$ & $181 $  \\ \hline
\hline
$6$ & $9$ & $x^6 + x^3 + 1$ & $19 $  \\ \hline
$6$ & $19$ & $x^6 + x^5 + 2 x^4 - 8 x^3 - x^2 + 5 x + 7$ & $7 $  \\ \hline
$6$ & $28$ & $x^6 + 5 x^4 + 6 x^2 + 1$ & $337 $  \\ \hline
$6$ & $39$ & $x^6 - x^5 + 5 x^4 + 6 x^3 + 15 x^2 + 4 x + 1$ & $31 $  \\ \hline
$6$ & $52$ & $x^6 + 13 x^4 + 26 x^2 + 13$ & $31 $  \\ \hline
$6$ & $63$ & $x^6 - 28 x^3 + 343$ & $193 $  \\ \hline
$6$ & $63$ & $x^6 + 35 x^3 + 343$ & $241 $  \\ \hline
$6$ & $91$ & $x^6 - x^5 + 21 x^4 - 22 x^3 + 58 x^2 + 23 x + 155$ & $31 $  \\ \hline
\hline
$8$ & $85$ & $x^8 - x^7 + 10 x^6 + 6 x^5 + 49 x^4 - 129 x^3 + 500 x^2 + 2044 x + 1616  $ & $433 $  \\ \hline
$8$ & $96$ & $x^8 + 24 x^6 + 180 x^4 + 432 x^2 + 162$ & $17 $  \\ \hline
\end{tabular}\\

\bigskip

All primes listed above have residue class degree $f=1$, but we have also found some examples of primes  (for characters of conductor $>100$) with $f=2$ (see below). Looking at the estimated probability (\ref{pred-prob}), it is not surprising that these primes are small.\\

\noindent \begin{tabular}{|l|l|l|l|} \hline
order & conductor & field (polynomial) & prime\\ \hline
$4$ & $187$ & $x^4 - x^3 + 45 x^2 + x + 562$ & $19 $ \\ \hline
\hline
$6$ & $259$ & $x^6 - x^5 + 22 x^4 + 398 x^3 + 1051 x^2 - x + 23605$ & $5 $  \\ \hline
\hline
$8$ & $187$ & $x^8 - x^7 + 44 x^6 - 45 x^5 + 423 x^4 - 877 x^3 + 1826 x^2 - 3515 x + 4591$ & $13 $  \\ \hline 
\end{tabular}\\

\bigskip

\subsection{Distribution of $\lambda$-invariants in the rank zero case} 

We provide data for the rank zero case and supplement the numerical data given in \cite{DeKn}. 
With the new formulas of this article, we computed the $\lambda$-invariants of $\chi = \theta \omega^ i$, where $\theta$ is odd of conductor $<100,000$ and $\theta(p) \neq 1$ if $i=1$. As in Section \ref{Sect41}, the first row gives the predicted distribution and the second row contains the number $N$ of tested characters as well as the computed proportions. \\

\noindent $p=3$, $\ord(\theta)=2$, twist $i=1$\\
\begin{tabular}{|l|l|l|l|l|l|l|l|l|} \hline
& $\lambda=0$ & $\lambda=1$ & $\lambda=2$ & $\lambda=3$ & $\lambda=4$ & $\lambda=5$ & $\lambda=6$ & $\lambda\geq 7$ \\ \hline
predicted & $0.5601$ & $0.2801$ & $0.1050$ & $0.0364$ & $0.0123$ & $0.0041$ & $0.0014$ & $0.0007$ \\ \hline
$N=18988$ & $0.6238$ & $0.2538$ & $0.0820$ & $0.0272$ & $0.0094$ & $0.0025$ & $0.0009$ & $0.0004$  \\ \hline
\end{tabular}\\

\noindent $p=3$, $\ord(\theta)=4$,  twist $i=1$\\
\begin{tabular}{|l|l|l|l|l|l|l|l|l|} \hline
& $\lambda=0$ & $\lambda=1$ & $\lambda=2$ & $\lambda=3$ & $\lambda=4$ & $\lambda=5$ & $\lambda=6$ & $\lambda\geq 7$ \\ \hline
predicted & $0.8766$ & $0.1096$ &  $0.0123$ & $0.0014$ &  $0.0002$ & $0.0000$ &  $0.0000$ &  $0.0000$  \\ \hline
$N=89374$ & $0.8880$ & $0.1004$ & $0.0106$ & $0.0009$ & $0.0002$ & $0.0000$ & $0.0000$ & $0.0000$ \\ \hline
\end{tabular}\\

\noindent $p=3$, $\ord(\theta)=8$,  twist $i=1$\\
\begin{tabular}{|l|l|l|l|l|l|l|l|l|} \hline
& $\lambda=0$ & $\lambda=1$ & $\lambda=2$ & $\lambda=3$ & $\lambda=4$ & $\lambda=5$ & $\lambda=6$ & $\lambda\geq 7$ \\ \hline
predicted & $0.8766$ & $0.1096$ &  $0.0123$ & $0.0014$ &  $0.0002$ & $0.0000$ &  $0.0000$ &  $0.0000$  \\ \hline
$N=129594$ &  $0.8819$ & $0.1052$ & $0.0114$ & $0.0013$ & $0.0002$ & $0.0000$ & $0.0000$ & $0.0000$ \\ \hline
\end{tabular}\\

\noindent $p=5$, $\ord(\theta)=2$,  twist $i=0$\\
\begin{tabular}{|l|l|l|l|l|l|l|l|l|} \hline
& $\lambda=0$ & $\lambda=1$ & $\lambda=2$ & $\lambda=3$ & $\lambda=4$ & $\lambda=5$ & $\lambda=6$ & $\lambda\geq 7$ \\ \hline
predicted & $0.7603$ &  $0.1901$ &  $0.0396$ &  $0.0080$ &  $0.0016$ & $0.0003$ &  $0.0001$ & $0.0000$ \\ \hline
$N=25324$ &
$0.7768$ & $0.1752$ & $0.0380$ & $0.0082$ & $0.0014$ & $0.0003$ & $0.0001$ & $0.0000$ \\ \hline
\end{tabular}\\

\noindent $p=5$, $\ord(\theta)=2$,  twist $i=1$\\
\begin{tabular}{|l|l|l|l|l|l|l|l|l|} \hline
& $\lambda=0$ & $\lambda=1$ & $\lambda=2$ & $\lambda=3$ & $\lambda=4$ & $\lambda=5$ & $\lambda=6$ & $\lambda\geq 7$ \\ \hline
predicted & $0.7603$ &  $0.1901$ &  $0.0396$ &  $0.0080$ &  $0.0016$ & $0.0003$ &  $0.0001$ & $0.0000$ \\ \hline
$N=17735$ & $0.7827$ & $0.1737$ & $0.0345$ & $0.0074$ & $0.0015$ & $0.0002$ & $0.0000$ & $0.0000$ \\ \hline
\end{tabular}\\

\noindent $p=5$, $\ord(\theta)=4$ , twist $i=1$\\
\begin{tabular}{|l|l|l|l|l|} \hline
$N$ & $\lambda=0$ & $\lambda=1$ & $\lambda=2$ & $\lambda\geq 3$ \\ \hline
predicted & $0.7603$ & $0.1901$ & $0.0396$ & $0.0100$ \\ \hline
$N=93424$ & $0.7760$ & $0.1790$ & $0.0360$ & $0.0090$ \\ \hline
\end{tabular}\\

\noindent $p=7$, $\ord(\theta)=3$, twist $i=0$\\
\begin{tabular}{|l|l|l|l|l|l|l|l|l|} \hline
& $\lambda=0$ & $\lambda=1$ & $\lambda=2$ & $\lambda=3$ & $\lambda=4$ & $\lambda=5$ & $\lambda=6$ & $\lambda\geq 7$ \\ \hline
predicted & $0.8368$ &  $0.1395$ & $0.0203$ &  $0.0029$ &  $0.0004$ &  $0.0001$ & $0.0000$ & $0.0000$  \\ \hline
$N=24676$ & $0.8407$ & $0.1332$ & $0.0221$ & $0.0035$ & $0.0004$ & $0.00004$ & $0.0000$ & $0.0000$ \\ \hline
\end{tabular}\\

\noindent $p=7$, $\ord(\theta)=4$,   twist $i=1$\\
\begin{tabular}{|l|l|l|l|l|} \hline
& $\lambda=0$ & $\lambda=1$ & $\lambda=2$ & $\lambda\geq 3$ \\ \hline
predicted & $0.9792$ & $0.0204$ & $0.0004$ & $0.0000$ \\ \hline
$N=85862$ & $0.9805$ & $0.0191$ & $0.0004$ & $0.0000$\\ \hline
\end{tabular}\\

\noindent $p=7$, $\ord(\theta)=6$,   twist $i=1$\\
\begin{tabular}{|l|l|l|l|l|} \hline
 & $\lambda=0$ & $\lambda=1$ & $\lambda=2$ & $\lambda\geq 3$ \\ \hline
predicted & $0.8368$ & $0.1395$ & $0.0203$ & $0.0034$  \\ \hline
$N=266777$ & $0.8430$ & $0.1343$ & $0.0197$ & $0.0031$
\\ \hline
\end{tabular}\\

\noindent $p=11$, $\ord(\theta)=6$,   twist $i=1$\\
\begin{tabular}{|l|l|l|l|l|} \hline
 & $\lambda=0$ & $\lambda=1$ & $\lambda=2$ & $\lambda\geq 3$ \\ \hline
predicted & $0.9917$ & $0.0083$ & $0.00007$ & $0.00000$ \\ \hline
$N=252356$ & $0.9924$ & $0.00756$ & $0.00004$ & $0.00000$\\ \hline
\end{tabular}\\

\backmatter

\bmhead{Acknowledgements}
I wish to thank Daniel Delbourgo for helpful discussions, his support and encouragement. I am also indebted to Chao Qin for his comments on an earlier version of the manuscript. Furthermore, I would like to thank the anonymous referees for their helpful suggestions.\\

\section*{Declarations}
\bmhead{Conflict of interest} The author has no competing interests to declare that are relevant to the content of this article.

\bmhead{Data availability}
All numerical data supporting the results of this study are presented in the tables within the article. The SageMath code used for the computations and to generate this data is publicly available on GitHub at \url{https://github.com/knospe/iwasawa}. 


\bibliography{special-values-p-adic-bibliography}


\begin{thebibliography}{16}
\ifx \bisbn   \undefined \def \bisbn  #1{ISBN #1}\fi
\ifx \binits  \undefined \def \binits#1{#1}\fi
\ifx \bauthor  \undefined \def \bauthor#1{#1}\fi
\ifx \batitle  \undefined \def \batitle#1{#1}\fi
\ifx \bjtitle  \undefined \def \bjtitle#1{#1}\fi
\ifx \bvolume  \undefined \def \bvolume#1{\textbf{#1}}\fi
\ifx \byear  \undefined \def \byear#1{#1}\fi
\ifx \bissue  \undefined \def \bissue#1{#1}\fi
\ifx \bfpage  \undefined \def \bfpage#1{#1}\fi
\ifx \blpage  \undefined \def \blpage #1{#1}\fi
\ifx \burl  \undefined \def \burl#1{\textsf{#1}}\fi
\ifx \doiurl  \undefined \def \doiurl#1{\url{https://doi.org/#1}}\fi
\ifx \betal  \undefined \def \betal{\textit{et al.}}\fi
\ifx \binstitute  \undefined \def \binstitute#1{#1}\fi
\ifx \binstitutionaled  \undefined \def \binstitutionaled#1{#1}\fi
\ifx \bctitle  \undefined \def \bctitle#1{#1}\fi
\ifx \beditor  \undefined \def \beditor#1{#1}\fi
\ifx \bpublisher  \undefined \def \bpublisher#1{#1}\fi
\ifx \bbtitle  \undefined \def \bbtitle#1{#1}\fi
\ifx \bedition  \undefined \def \bedition#1{#1}\fi
\ifx \bseriesno  \undefined \def \bseriesno#1{#1}\fi
\ifx \blocation  \undefined \def \blocation#1{#1}\fi
\ifx \bsertitle  \undefined \def \bsertitle#1{#1}\fi
\ifx \bsnm \undefined \def \bsnm#1{#1}\fi
\ifx \bsuffix \undefined \def \bsuffix#1{#1}\fi
\ifx \bparticle \undefined \def \bparticle#1{#1}\fi
\ifx \barticle \undefined \def \barticle#1{#1}\fi
\bibcommenthead
\ifx \bconfdate \undefined \def \bconfdate #1{#1}\fi
\ifx \botherref \undefined \def \botherref #1{#1}\fi
\ifx \url \undefined \def \url#1{\textsf{#1}}\fi
\ifx \bchapter \undefined \def \bchapter#1{#1}\fi
\ifx \bbook \undefined \def \bbook#1{#1}\fi
\ifx \bcomment \undefined \def \bcomment#1{#1}\fi
\ifx \oauthor \undefined \def \oauthor#1{#1}\fi
\ifx \citeauthoryear \undefined \def \citeauthoryear#1{#1}\fi
\ifx \endbibitem  \undefined \def \endbibitem {}\fi
\ifx \bconflocation  \undefined \def \bconflocation#1{#1}\fi
\ifx \arxivurl  \undefined \def \arxivurl#1{\textsf{#1}}\fi
\csname PreBibitemsHook\endcsname

\bibitem[\protect\citeauthoryear{Kubota and Leopoldt}{1964}]{KuLe}
\begin{barticle}
\bauthor{\bsnm{Kubota}, \binits{T.}},
\bauthor{\bsnm{Leopoldt}, \binits{H.-W.}}:
\batitle{{Eine $p$-adische Theorie der Zetawerte, I: Einf\"{u}hrung der
  $p$-adischen Dirichletschen L-Funktionen}}.
\bjtitle{J. Reine Angew. Math.}
\bvolume{214},
\bfpage{328}--\blpage{339}
(\byear{1964})
\end{barticle}
\endbibitem

\bibitem[\protect\citeauthoryear{Ferrero and Washington}{1979}]{FeWa}
\begin{barticle}
\bauthor{\bsnm{Ferrero}, \binits{B.}},
\bauthor{\bsnm{Washington}, \binits{L.C.}}:
\batitle{{The Iwasawa $\mu_p$-invariant vanishes for abelian number fields}}.
\bjtitle{Annals of Math.}
\bvolume{109},
\bfpage{377}--\blpage{395}
(\byear{1979})
\end{barticle}
\endbibitem

\bibitem[\protect\citeauthoryear{Delbourgo and Knospe}{2023}]{DeKn}
\begin{barticle}
\bauthor{\bsnm{Delbourgo}, \binits{D.}},
\bauthor{\bsnm{Knospe}, \binits{H.}}:
\batitle{{On Iwasawa $\lambda$-invariants for abelian number fields and random
  matrix heuristics}}.
\bjtitle{Math. Comp.}
\bvolume{92},
\bfpage{1817}--\blpage{1836}
(\byear{2023})
\end{barticle}
\endbibitem

\bibitem[\protect\citeauthoryear{Ellenberg et~al.}{2011}]{EJV}
\begin{barticle}
\bauthor{\bsnm{Ellenberg}, \binits{J.}},
\bauthor{\bsnm{Jain}, \binits{S.}},
\bauthor{\bsnm{Venkatesh}, \binits{A.}}:
\batitle{{Modelling $\lambda$-invariants by $p$-adic random matrices}}.
\bjtitle{Commun. Pure Appl. Math.}
\bvolume{64},
\bfpage{1243}--\blpage{1262}
(\byear{2011})
\end{barticle}
\endbibitem

\bibitem[\protect\citeauthoryear{Ferrero and Greenberg}{1978}]{FeGr}
\begin{barticle}
\bauthor{\bsnm{Ferrero}, \binits{B.}},
\bauthor{\bsnm{Greenberg}, \binits{R.}}:
\batitle{{On the Behavior of $p$-adic $L$-functions at $s=0$}}.
\bjtitle{Inv. Math.}
\bvolume{50},
\bfpage{91}--\blpage{102}
(\byear{1978})
\end{barticle}
\endbibitem

\bibitem[\protect\citeauthoryear{Gross and Dasgupta}{2023}]{GrDa}
\begin{botherref}
\oauthor{\bsnm{Gross}, \binits{B.H.}},
\oauthor{\bsnm{Dasgupta}, \binits{S.}}:
{Two encounters with the $p$-adic Stark conjecture}.
Preprint at \url{https://arxiv.org/abs/2303.03299}
(2023)
\end{botherref}
\endbibitem

\bibitem[\protect\citeauthoryear{Gross and Koblitz}{1979}]{GrKo}
\begin{barticle}
\bauthor{\bsnm{Gross}, \binits{B.H.}},
\bauthor{\bsnm{Koblitz}, \binits{N.}}:
\batitle{{Gauss sums and the $p$-adic Gamma function}}.
\bjtitle{Annals of Math.}
\bvolume{109},
\bfpage{569}--\blpage{581}
(\byear{1979})
\end{barticle}
\endbibitem

\bibitem[\protect\citeauthoryear{Washington}{1981}]{Wa0}
\begin{barticle}
\bauthor{\bsnm{Washington}, \binits{L.C.}}:
\batitle{{The derivative of $p$-adic $L$-functions}}.
\bjtitle{Acta Arithmetica}
\bvolume{XL},
\bfpage{109}--\blpage{115}
(\byear{1981})
\end{barticle}
\endbibitem

\bibitem[\protect\citeauthoryear{Lang}{1990}]{lang}
\begin{bbook}
\bauthor{\bsnm{Lang}, \binits{S.}}:
\bbtitle{Cyclotomic Fields I and II},
\bedition{Combined second} edn.
\bsertitle{Graduate Texts in Mathematics},
vol. \bseriesno{121}
(\byear{1990})
\end{bbook}
\endbibitem

\bibitem[\protect\citeauthoryear{Diamond}{1977}]{Dia}
\begin{barticle}
\bauthor{\bsnm{Diamond}, \binits{J.}}:
\batitle{{The $p$-adic log gamma function and $p$-adic Euler constants}}.
\bjtitle{Trans. Amer. Math. Soc.}
\bvolume{233},
\bfpage{321}--\blpage{337}
(\byear{1977})
\end{barticle}
\endbibitem

\bibitem[\protect\citeauthoryear{Washington}{1996}]{Wa2}
\begin{bbook}
\bauthor{\bsnm{Washington}, \binits{L.C.}}:
\bbtitle{Introduction to Cyclotomic Fields},
\bedition{2}nd edn.
\bsertitle{Graduate Texts in Mathematics},
vol. \bseriesno{83}
(\byear{1996})
\end{bbook}
\endbibitem

\bibitem[\protect\citeauthoryear{Stokes}{2023}]{Sto}
\begin{barticle}
\bauthor{\bsnm{Stokes}, \binits{M.}}:
\batitle{{On Gauss factorials and their connection to the cyclotomic
  $\lambda$-invariants of imaginary quadratic fields}}.
\bjtitle{Journal of Number Theory}
\bvolume{246},
\bfpage{279}--\blpage{293}
(\byear{2023})
\end{barticle}
\endbibitem

\bibitem[\protect\citeauthoryear{Greenberg}{2001}]{Gr}
\begin{bchapter}
\bauthor{\bsnm{Greenberg}, \binits{R.}}:
\bctitle{{Iwasawa theory - Past and Present}}.
In: \bbtitle{Class Field Theory - Its Centenary and Prospect}
vol. \bseriesno{30},
pp. \bfpage{335}--\blpage{386}.
\bpublisher{Mathematical Society of Japan},
\blocation{Tokyo}
(\byear{2001})
\end{bchapter}
\endbibitem

\bibitem[\protect\citeauthoryear{Dummit et~al.}{1991}]{Du}
\begin{barticle}
\bauthor{\bsnm{Dummit}, \binits{D.S.}},
\bauthor{\bsnm{Ford}, \binits{D.}},
\bauthor{\bsnm{Kisilevsky}, \binits{H.}},
\bauthor{\bsnm{Sands}, \binits{J.W.}}:
\batitle{{Computation of Iwasawa Lambda Invariants for Imaginary Quadratic
  Fields}}.
\bjtitle{Journal of Number Theory}
\bvolume{37},
\bfpage{100}--\blpage{121}
(\byear{1991})
\end{barticle}
\endbibitem

\bibitem[\protect\citeauthoryear{Ernvall and Metsänkylä}{1987}]{ErMe}
\begin{barticle}
\bauthor{\bsnm{Ernvall}, \binits{R.}},
\bauthor{\bsnm{Metsänkylä}, \binits{T.}}:
\batitle{{A Method for Computing the Iwasawa $\lambda$-Invariant}}.
\bjtitle{Math. Comp.}
\bvolume{49},
\bfpage{281}--\blpage{294}
(\byear{1987})
\end{barticle}
\endbibitem

\bibitem[\protect\citeauthoryear{Knospe and Washington}{2021}]{KnWa}
\begin{barticle}
\bauthor{\bsnm{Knospe}, \binits{H.}},
\bauthor{\bsnm{Washington}, \binits{L.C.}}:
\batitle{{Dirichlet series expansions of $p$-adic $L$-functions}}.
\bjtitle{Abhandlungen aus dem Mathematischen Seminar der Universit\"{a}t
  Hamburg}
\bvolume{91},
\bfpage{325}--\blpage{334}
(\byear{2021})
\end{barticle}
\endbibitem

\end{thebibliography}

\end{document}